# BINOMIAL TRANSFORM OF SEQUENCES COUNTING $N$-ARY CONVEXITIES

Aidar Dulliev[*], Daniil Naumikhin[*]

June 8, 2026

## Abstract

We consider the enumeration of $N$-ary convex structures on finite sets. Our main result shows that the total number of convexities is expressed as the binomial transform of the sequence of numbers of grounded convexities. We present the exact numbers of all $N$-ary and grounded $N$-ary convexities for $|X| \leq 5$. The obtained integer sequences have been cross-referenced with the OEIS. As a result, we identified both known sequences (such as A000798) and new ones not previously listed.

Total: 19 pages (translated English text: pp. 1–13; original Russian text: pp. 14–19).

[*] Department of Applied Mathematics and Informatics, Kazan National Research Technical University named after A.N. Tupolev-KAI, Russian Federation. E-mail: {dulliev@yandex.ru, dan.naumihin@yandex.ru}.

We sincerely apologize for the English text below. According to the changed policy, we urgently need to translate the text from Russian using a machine translator. The original text in Russian is presented after the English version. Both versions – the English and Russian – are planned to be corrected further.

# 1. INTRODUCTION

Let $X$ be an arbitrary set. A *convexity* $\mathbf{G}$ in $X$ is a family of subsets of $X$ closed under arbitrary intersections of its elements and containing $X$; in this case, the pair $(X, \mathbf{G})$ is called *a convexity space (a convex structure)*. Every convexity is in one-to-one correspondence with *a convex hull operator* $g: \mathcal{P}(X) \to \mathcal{P}(X)$ satisfying, for all $A \subset B \subset X$ the axioms: $A \subset gA$, $gA \subset gB$, $ggA = gA$; the correspondence is established by the rules: $gA = \bigcap_{A \subset B \in \mathbf{G}} B$ and $\mathbf{G} = \{A \subset X : gA = A\}$. The convexity $\mathbf{G}$ is called $N$-ary if $A \in \mathbf{G} \Leftrightarrow \forall B \subset A : (|B| \leq N) \Rightarrow gB \subset A$ whenever $A \subset X$. We shall say that a convexity $\mathbf{G}$ satisfies the grounding axiom (such a convexity will also be called a convexity with grounding, or a grounded convexity) if the following condition holds $\emptyset \in \mathbf{G}$.

Two monographs, now classical, are devoted to a systematic exposition of the foundations of convexity theory: V.P. Soltan [1] and M.L.J. van de Vel [2]). It should be noted that the definitions of convex spaces in them are different (non-equivalent): in [1] a more general definition is given, while in [2] the inductivity condition for convexity is added. Questions related to the enumeration of convex structures of various types and of structures generated by convex ones on finite carriers have been studied very little. The results obtained to date mainly reduce to the enumeration, not complete, of other known structures for which equivalence to certain convex structures has been proved and enumeration has been carried out. The collection of numbers, each of which is equal to the number of convex structures of one kind or another on a finite carrier and depends on the number of its elements, naturally forms an integer sequence. Sometimes it is possible to find a closed formula for an arbitrary term of such a sequence. In the present paper a relation is obtained that connects the total number of convexities with the number of grounded convexities; it is expressed as a binomial transform of the corresponding sequences.

The integer sequences obtained from the developed computer program are compared with sequences from the international “Online Encyclopedia of Integer Sequences”. Several sequences, as expected, turned out to be new.

We give the list of symbols and notation used in the paper:

- $|A|$ – the cardinality of the set $A$;
- $\mathcal{P}(X)$ – the family of all subsets of the set $X$, the power set of $X$;
- $\mathbf{G}$ – a convexity in $X$;
- $gA = \cap \{B \in \mathbf{G} : A \subset B\}$ – the convex hull of the set $A \subset X$
- $\Gamma^N(X)$ – the family of all $N$-ary convexities on the carrier $X$;
- $\Gamma_0^N(X)$ – the family of all grounded $N$-ary convexities on the carrier $X$;
- $X_n$ – finite set of cardinality *n*.

# 2. PROPERTIES OF $\Gamma^N(X)$, $\Gamma_0^N(X)$

Everywhere below, it is assumed that *X* is a finite set.

**Lemma 1.** *Let* $\mathbf{G} \in \Gamma^N(X)$, $C$ *is the smallest convex set*

$$C = g\emptyset = \bigcap_{A \in \mathbf{G}} A,$$

*and set* $\mathbf{H}$ *is constructed* according to *the rule*

$$\mathbf{H} = \{A \backslash C : A \in \mathbf{G}\}. \quad (1)$$

*Then between* $\mathbf{G}$ *and* $\mathbf{H}$ *there exists a one-to-one correspondence.*

**Proof.** It is clear from the construction that the mapping $\varphi: X \mapsto X\backslash C$ is surjective:

$$\forall B \in \mathbf{H} \exists A \in \mathbf{G}: B = A\backslash C.$$

We prove injectivity of $\varphi$. Suppose the contrary: there are two distinct sets $A_1, A_2$, for which $\varphi A_1 = \varphi A_2$, i.e.

$$A_1 \neq A_2, A_1\backslash C = A_2\backslash C. \tag{2}$$

From the construction of $C$:

$$\forall A \in \mathbf{G}: C \subset A. \tag{3}$$

Consider an arbitrary point $x$ from $A_1$. Suppose $x$ belongs also to $C$, then by (3) $x$ also belongs to $A_2$. Otherwise, if $x \in A_1\backslash C$, then by (2) $x$ also belongs to $A_2\backslash C$. Then $A_1 \subset A_2$. Repeating the same reasoning for $A_2$, we obtain $A_2 \subset A_1$. Thus, $A_1 = A_2$, which contradicts construction (2).

The mapping $X \mapsto X\backslash C$ is both injective and surjective, and hence bijective. ■

**Lemma 2.** *Let* $\mathbf{G} \in \Gamma^N(X)$, *then* $\mathbf{H}$, *constructed by rule (1) is a grounded convexity on* $Y = X\backslash C$.

**Proof.** 1. It is obvious that according to (1) we have: $Y \in \mathbf{H}$.

2. This family of sets is closed under intersection

$$\forall B_1, B_2 \in \mathbf{H} \exists A_1, A_2 \in \mathbf{G}: B_1 = A_1\backslash C, B_2 = A_2\backslash C;$$
$$A_1 \cap A_2 \in \mathbf{G} \Rightarrow (A_1 \cap A_2)\backslash C = B_1 \cap B_2 \in \mathbf{H}.$$

3. This family of sets contains the empty set by construction:

$$g\emptyset \in G \Rightarrow g\emptyset\backslash g\emptyset = \emptyset \in \mathbf{H}.$$ ■

**Lemma 3.** *Let* $\mathbf{G} \in \Gamma^N(X)$. *The convex hull* $h$ *of the convexity* $\mathbf{H}$, *constructed by rule (1), has the following relation to the convex hull* $g$:

$$\forall A \subset Y: gA = C \cup hA \tag{4}$$

*or*

$$hA = gA\backslash C. \tag{5}$$

*Moreover, the following representation holds*

$$A \subset X \Rightarrow gA = h(A\backslash C) \cup C. \tag{6}$$

**Proof.** The validity of relation (4) follows from the following chain of equalities:

$$gA = \bigcap_{\substack{B \in \mathbf{G} \\ A \subset B}} B = \bigcap_{\substack{Z \cup C \in \mathbf{G} \\ A \subset Z \cup C}} (Z \cup C) =$$
$$= C \cup \bigcap_{\substack{Z \cup C \in \mathbf{G} \\ A \subset Z \cup C}} Z = C \cup \bigcap_{\substack{Z \in \mathbf{H} \\ A \subset Z}} Z = C \cup hA.$$

By construction of the convex hull $h$ we have:

$$hA \in \mathbf{H}, \mathbf{H} = \{B\backslash C: B \in \mathbf{G}\}.$$

hence $hA \cap C = \emptyset$, and therefore, from what has already been proved above, relation (5).

To prove relation (6) consider an arbitrary $A \subset X$. Note that $gA = g(A\backslash C)$. Indeed, adding to $A$ elements from $C$ cannot enlarge the convex hull, since $C = g\emptyset \subset g(A\backslash C)$ (the minimal convex set $C$ is contained in every convex hull). Therefore, for every $x \in A \cap C$ we have $x \in C \subset g(A\backslash C)$, and hence $gA = g(A\backslash C)$.

Now $A\backslash C \subset Y = X\backslash C$, so formula (4) can be applied to the set $A\backslash C$:

$$g(A\backslash C) = C \cup h(A\backslash C).$$

Thus,

$$gA = C \cup h(A\backslash C),$$

which exactly corresponds to formula (6). ■

**Lemma 4.** *Let* $\mathbf{G} \in \Gamma^N(X)$ . *The convexity* $\mathbf{H}$*, constructed by rule (1), is* $N$*-ary.*

**Proof.** Consider an arbitrary convexity $\mathbf{G}$ from $\Gamma^N(X)$. Then, by the definition of a $N$-ary convexity, we have:

$$\forall A \subset X: [A \in \mathbf{G} \Leftrightarrow (\forall B \subset A: |B| \le N \Rightarrow gB \subset A)].$$

From Lemma 2 and the relation between the convex hulls $g$ and $h$ (Lemma 3), restricting the property under consideration to the set $Y = X\backslash C$:

$$\forall A \subset Y: [A \cup C \in \mathbf{G} \Leftrightarrow (\forall B \subset C \cup A: |B| \le N \Rightarrow hB \cup C \subset A \cup C)].$$

Applying the mapping $\mathbf{G}$ by $\mathbf{H}$, and also taking into account the construction of the convex hull $h$, namely the fact that $\mathbf{H}$ is constructed as a subset of $X\backslash C$, it follows that every convex hull $h$ is contained entirely in $X\backslash C$, exactly as the set $A$, belonging to $\mathbf{H}$ on the left-hand side; that is, the following holds:

$$\forall A \subset Y[A \in \mathbf{H} \Leftrightarrow (\forall B \subset A: |B| \le N \Rightarrow hB \subset A)].$$

This corresponds to the definition of a $N$-ary set on $Y = X\backslash C$. From Lemma 2 proved above, this convexity is also grounded, so $\mathbf{H}$ is indeed a grounded convexity on the carrier $Y = X\backslash C$. ■

Let $C \subset X$ be arbitrary. Define by $E_C$ the family of convexities from $\Gamma^N(X)$ as follows:

$$\mathbf{G} \in E_C \Leftrightarrow g\emptyset = C.$$

The following lemma holds.

**Lemma 5.** *To every convexity* $\mathbf{G} \in E_C$ *one can assign a unique convexity* $\boldsymbol{H} \in \Gamma_0^N(X\backslash C)$*, defined by rule (1).*

*Conversely, for every* $\mathbf{H} \in \Gamma_0^N(X\backslash C)$ *the family*

$$\mathbf{G} = \{C \cup D: D \in \mathbf{H}\}$$

*belongs to* $E_C$.

*The constructed correspondence establishes between* $E_C$ *and* $\Gamma_0^N(X\backslash C)$ *a bijection* $\xi_C$.

**Proof.** (*direct mapping*): Let $\mathbf{G} \in E_C$. By assumption $g\emptyset = C$. Construct $\mathbf{H}$ according to rule (1). Lemmas 2 and 3 guarantee that $\mathbf{H} \in \Gamma_0^N(X\backslash C)$. Lemma 1 guarantees that the constructed mapping is a bijection.

(*inverse mapping*): Let $\mathbf{H} \in \Gamma_0^N(Y)$, where $Y = X\backslash C$. Define $\mathbf{G} = \{C \cup D: D \in \mathbf{H}\}$. Let us verify that $\mathbf{G} \in E_C$:

1. By construction $C = C \cup \emptyset \in \mathbf{G}$, since $\emptyset \in \mathbf{H}$. Moreover, for every $C \cup D \in \mathbf{G}$ we have $C \subset C \cup D$. Therefore, $C$ is contained in all elements of $\mathbf{G}$, and hence $C \subset g\emptyset$.
2. Let us show that $g\emptyset \subset C$. Suppose $x \in g\emptyset$. Then $x$ belongs to all elements of $\mathbf{G}$, including $C$. Therefore $x \in C$.
3. From the previous points $g\emptyset = C$.
4. Let us verify that $\mathbf{G}$ is an $N$-ary convexity:

(a) $\mathbf{G}$ is closed under intersections. For any $C \cup D_1, C \cup D_2 \in \mathbf{G}$:

$$(C \cup D_1) \cap (C \cup D_2) = C \cup (D_1 \cap D_2).$$

Since $D_1, D_2 \in \mathbf{H}$, $D_1 \cap D_2 \in \mathbf{H}$. Therefore $C \cup (D_1 \cap D_2) \in \mathbf{G}$.

(b) The carrier of the convexity is $X$. Since $Y \in \mathbf{H}$, $X = C \cup Y \in \mathbf{G}$ by construction. The minimal convex element of the convexity $\mathbf{G}$ is $C$, since $\emptyset \in \mathbf{H}$, and by construction $C = C \cup \emptyset$.

(c) $N$-arity. Take an arbitrary $A = C \cup D \in \mathbf{G}(D \in \mathbf{H})$. Suppose $B \subset A, |B| \le N$. Split $B = B_C \cup B_D, B_C \subset C, B_D \subset D$. The convex hull $\mathbf{G}$ has the following property (a consequence of

Lemma 3 applied to $\mathbf{H}$ and the inverse construction):

$$gB = g(B_C \cup B_D) = C \cup hB_D.$$

Since $C$ is the minimal convex set, $B_D \subset Y$, then $gB_D = C \cup hB_d$. Since $|B_D| \leq N, D \in \mathbf{H}$ — $N$-ary, then $hB_D \subset D$. Hence

$$gB = C \cup hB_D \subset C \cup D = A.$$

Therefore, $A$ satisfies the definition of a $N$-ary set, and therefore $A \in \mathbf{G}$.

Conversely, if $A = C \cup D$ satisfies the $N$-ary condition, then for every $B_D \subset D$ with $|B_D| \leq N$ we have $gB_D \subset A$, whence $C \cup hB_D \subset C \cup D$, i.e. $hB_D \subset D$. Hence $D$ satisfies the condition of $N$-arity in $\mathbf{H}$, since $\mathbf{H}$ is $N$-ary convexity, $D \in \mathbf{H}$. Thus, $\mathbf{G}$ consists exactly of the sets $C \cup D$, $D \in \mathbf{H}$.

(d) Uniqueness. If $\mathbf{G}$ is an arbitrary $N$-ary convexity with $g\emptyset = C$, $\mathbf{H} = \{A \setminus C : A \in \mathbf{G}\}$, then by construction $\mathbf{G}$ must coincide with $\{C \cup D : D \in \mathbf{H}\}$. Hence, the correspondence between $\mathbf{H}$ and $\mathbf{G}$ is one-to-one.

Thus, the constructed correspondence establishes a bijection between $E_C$ and $\Gamma_0^N(X \setminus C)$. ■

**Proposition.** *The total number of $N$-ary convex structures is expressed via the binomial transform of the number of grounded $N$-ary convex structures. Namely,*

$$|\Gamma^N(X_n)| = \sum_{k=0}^{n} C_n^k |\Gamma_0^N(X_k)|.$$

**Proof.** The classes $E_C$ for distinct $C \subset X$ form a partition of the set $\Gamma^N(X)$. The cardinality of each class is $|\Gamma_0^N(X \setminus C)|$, therefore

$$|\Gamma^N(X)| = \sum_{C \subset X} |\Gamma_0^N(X \setminus C)|. \tag{7}$$

The case $C = \emptyset$ in formula (7) corresponds to one summand $|\Gamma_0^N(X_n)|$. The case $C = \{x\}, x \in X_n$ in formula (7) corresponds to $n$ such summands $|\Gamma_0^N(X_{n-1})|$. Similarly, the case $|C| = k, 2 \leq k \leq n$, in formula (7) corresponds to $C_n^k$ such summands $|\Gamma_0^N(X_{n-k})|$ (one can choose $C_n^k$ different subsets $C \subset X$ of cardinality $k$). Thus,

$$\sum_{C \subset X} |\Gamma_0^N(X \setminus C)| = \sum_{k=0}^{n} C_n^k |\Gamma_0^N(X_{n-k})|.$$

Replacing $k$ by $n - k$ and using the property $C_n^k = C_n^{n-k}$, we obtain the required result

$$|\Gamma^N(X)| = \sum_{k=0}^{n} C_n^{n-k} |\Gamma_0^N(X_k)| = \sum_{k=0}^{n} C_n^k |\Gamma_0^N(X_k)|.$$ ■

## 3. NUMERICAL EXAMPLES

On the basis of the developed algorithm and program, a complete enumeration was carried out of $N$-ary convexities on finite carriers of cardinality $n = 0,1,\ldots,6$ for arities $N = 0,1,\ldots,6$. The results of the numerical experiment are presented in Tables 1 and 2. For $n = 0$ the unique empty set is meant, and for $n = 1$ the one-element carrier is meant. For all investigated arities and dimensions, therelation holds, established in Proposition in the second section.

**Table 1. Number of $N$-ary convexities, $|\Gamma^N(X_n)|$**

| $N$ \ $n$ | 0 | 1 | 2 | 3 | 4 | 5 | 6 | OEIS |
|---|---|---|---|---|---|---|---|---|
| **0** | 1 | 2 | 4 | 8 | 16 | 32 | 64 | A000079 |
| **1** | 1 | 2 | 7 | 45 | 500 | 9053 | 257151 | A326878 |
| **2** | 1 | 2 | 7 | 61 | 2271 | 600408 | 1556743050 | - |
| **3** | 1 | 2 | 7 | 61 | 2480 | 1373701 | - | - |
| **4** | 1 | 2 | 7 | 61 | 2480 | 1385552 | - | - |

**Table 2. Number of grounded $N$-ary convexities, $|\Gamma_0^N(X_n)|$**

| $N$ \ $n$ | 0 | 1 | 2 | 3 | 4 | 5 | 6 | OEIS |
|---|---|---|---|---|---|---|---|---|
| **0** | 1 | 1 | 1 | 1 | 1 | 1 | 1 | A000012 |
| **1** | 1 | 1 | 4 | 29 | 355 | 6942 | 209527 | A000798 |
| **2** | 1 | 1 | 4 | 45 | 2062 | 589602 | 1553173541 | A364656 |
| **3** | 1 | 1 | 4 | 45 | 2271 | 1361850 | - | A395658 |
| **4** | 1 | 1 | 4 | 45 | 2271 | 1373701 | - | - |

In Table 3 all *binary* convexities constructed on carriers with 0, 1, 2, and 3 elements are given. To the right of each of them after the vertical bar | is indicated the corresponding groundingconstructed by excluding the minimal convex set $g\emptyset$ according to rule (1). The empty set everywhere is denoted by the symbol {}. In Table 4 gives all *ternary,* but *not binary* convexities constructed on a carrier with 4 elements.

**Table 3. Binary convexities**

| **$n = 0$** | |
|---|---|
| **No** | **G \| H** |
| 1 | {{}} \| {{}} |

| **$n = 1$** | |
|---|---|
| **No** | **G \| H** |
| 1 | {{0}} \| {{}} |
| 2 | {{0}, {}} \| {{0}, {}} |

| **$n = 2$** | |
|---|---|
| **No** | **G \| H** |
| 1 | {{0, 1}} \| {{}} |
| 2 | {{0, 1}, {0}} \| {{1}, {}} |
| 3 | {{0, 1}, {}, {0}} \| {{0, 1}, {}, {0}} |
| 4 | {{0, 1}, {}} \| {{0, 1}, {}} |
| 5 | {{0, 1}, {1}} \| {{0}, {}} |
| 6 | {{0, 1}, {}, {1}} \| {{0, 1}, {}, {1}} |
| 7 | {{0, 1}, {}, {0}, {1}} \| {{0, 1}, {}, {0}, {1}} |

<table>
<tr><th colspan="2">$n = 3$</th></tr>
<tr><th>No</th><th>G | H</th></tr>
<tr><td>1</td><td>{{0, 1, 2}, {}, {0, 2}} | {{0, 1, 2}, {}, {0, 2}}</td></tr>
<tr><td>2</td><td>{{0, 1, 2}, {0}, {0, 2}} | {{1, 2}, {}, {2}}</td></tr>
<tr><td>3</td><td>{{0, 1, 2}, {}, {0}, {0, 2}} | {{0, 1, 2}, {}, {0}, {0, 2}}</td></tr>
<tr><td>4</td><td>{{0, 1, 2}, {}, {0}, {1}, {2}, {0, 2}} | {{0, 1, 2}, {}, {0}, {1}, {2}, {0, 2}}</td></tr>
<tr><td>5</td><td>{{0, 1, 2}, {}, {0}, {0, 1}, {2}, {0, 2}} | {{0, 1, 2}, {}, {0}, {0, 1}, {2}, {0, 2}}</td></tr>
<tr><td>6</td><td>{{0, 1, 2}, {}, {0}, {1}, {0, 1}, {2}, {0, 2}} | {{0, 1, 2}, {}, {0}, {1}, {0, 1}, {2}, {0, 2}}</td></tr>
<tr><td>7</td><td>{{0, 1, 2}, {1, 2}} | {{0}, {}}</td></tr>
<tr><td>8</td><td>{{0, 1, 2}, {}, {1, 2}} | {{0, 1, 2}, {}, {1, 2}}</td></tr>
<tr><td>9</td><td>{{0, 1, 2}} | {{}}</td></tr>
<tr><td>10</td><td>{{0, 1, 2}, {}} | {{0, 1, 2}, {}}</td></tr>
<tr><td>11</td><td>{{0, 1, 2}, {}, {1}, {0, 1}, {2}, {1, 2}} | {{0, 1, 2}, {}, {1}, {0, 1}, {2}, {1, 2}}</td></tr>
<tr><td>12</td><td>{{0, 1, 2}, {}, {0}, {1}, {0, 1}, {2}, {1, 2}} | {{0, 1, 2}, {}, {0}, {1}, {0, 1}, {2}, {1, 2}}</td></tr>
<tr><td>13</td><td>{{0, 1, 2}, {}, {0}, {0, 1}} | {{0, 1, 2}, {}, {0}, {0, 1}}</td></tr>
<tr><td>14</td><td>{{0, 1, 2}, {}, {0}, {1, 2}} | {{0, 1, 2}, {}, {0}, {1, 2}}</td></tr>
<tr><td>15</td><td>{{0, 1, 2}, {1}, {1, 2}} | {{0, 2}, {}, {2}}</td></tr>
<tr><td>16</td><td>{{0, 1, 2}, {}, {1}, {1, 2}} | {{0, 1, 2}, {}, {1}, {1, 2}}</td></tr>
<tr><td>17</td><td>{{0, 1, 2}, {}, {0}, {1}, {1, 2}} | {{0, 1, 2}, {}, {0}, {1}, {1, 2}}</td></tr>
<tr><td>18</td><td>{{0, 1, 2}, {1}, {0, 1}, {1, 2}} | {{0, 2}, {}, {0}, {2}}</td></tr>
<tr><td>19</td><td>{{0, 1, 2}, {}, {1}, {0, 1}, {1, 2}} | {{0, 1, 2}, {}, {1}, {0, 1}, {1, 2}}</td></tr>
<tr><td>20</td><td>{{0, 1, 2}, {}, {0}, {1}, {0, 1}, {1, 2}} | {{0, 1, 2}, {}, {0}, {1}, {0, 1}, {1, 2}}</td></tr>
<tr><td>21</td><td>{{0, 1, 2}, {2}, {1, 2}} | {{0, 1}, {}, {1}}</td></tr>
<tr><td>22</td><td>{{0, 1, 2}, {}, {2}, {1, 2}} | {{0, 1, 2}, {}, {2}, {1, 2}}</td></tr>
<tr><td>23</td><td>{{0, 1, 2}, {}, {0}, {2}, {1, 2}} | {{0, 1, 2}, {}, {0}, {2}, {1, 2}}</td></tr>
<tr><td>24</td><td>{{0, 1, 2}, {}, {1}, {2}, {1, 2}} | {{0, 1, 2}, {}, {1}, {2}, {1, 2}}</td></tr>
<tr><td>25</td><td>{{0, 1, 2}, {}, {0}, {1}, {2}, {1, 2}} | {{0, 1, 2}, {}, {0}, {1}, {2}, {1, 2}}</td></tr>
<tr><td>26</td><td>{{0, 1, 2}, {1}, {0, 1}} | {{0, 2}, {}, {0}}</td></tr>
<tr><td>27</td><td>{{0, 1, 2}, {}, {1}, {0, 1}} | {{0, 1, 2}, {}, {1}, {0, 1}}</td></tr>
<tr><td>28</td><td>{{0, 1, 2}, {}, {0}, {1}, {2}} | {{0, 1, 2}, {}, {0}, {1}, {2}}</td></tr>
<tr><td>29</td><td>{{0, 1, 2}, {}, {0, 1}, {2}} | {{0, 1, 2}, {}, {0, 1}, {2}}</td></tr>
<tr><td>30</td><td>{{0, 1, 2}, {}, {0}, {0, 1}, {2}} | {{0, 1, 2}, {}, {0}, {0, 1}, {2}}</td></tr>
<tr><td>31</td><td>{{0, 1, 2}, {}, {1}, {0, 1}, {2}} | {{0, 1, 2}, {}, {1}, {0, 1}, {2}}</td></tr>
<tr><td>32</td><td>{{0, 1, 2}, {}, {0}, {1}, {0, 1}, {2}} | {{0, 1, 2}, {}, {0}, {1}, {0, 1}, {2}}</td></tr>
<tr><td>33</td><td>{{0, 1, 2}, {0, 2}} | {{1}, {}}</td></tr>
<tr><td>34</td><td>{{0, 1, 2}, {}, {0}, {1}, {2}, {0, 2}, {1, 2}} | {{0, 1, 2}, {}, {0}, {1}, {2}, {0, 2}, {1, 2}}</td></tr>
<tr><td>35</td><td>{{0, 1, 2}, {}, {0}, {1}, {0, 1}, {2}, {0, 2}, {1, 2}} | {{0, 1, 2}, {}, {0}, {1}, {0, 1}, {2}, {0, 2}, {1, 2}}</td></tr>
<tr><td>36</td><td>{{0, 1, 2}, {0}} | {{1, 2}, {}}</td></tr>
<tr><td>37</td><td>{{0, 1, 2}, {}, {0}} | {{0, 1, 2}, {}, {0}}</td></tr>
<tr><td>38</td><td>{{0, 1, 2}, {1}} | {{0, 2}, {}}</td></tr>
<tr><td>39</td><td>{{0, 1, 2}, {}, {1}} | {{0, 1, 2}, {}, {1}}</td></tr>
<tr><td>40</td><td>{{0, 1, 2}, {}, {0}, {1}, {0, 1}, {0, 2}} | {{0, 1, 2}, {}, {0}, {1}, {0, 1}, {0, 2}}</td></tr>
<tr><td>41</td><td>{{0, 1, 2}, {2}, {0, 2}} | {{0, 1}, {}, {0}}</td></tr>
<tr><td>42</td><td>{{0, 1, 2}, {2}, {0, 2}, {1, 2}} | {{0, 1}, {}, {0}, {1}}</td></tr>
<tr><td>43</td><td>{{0, 1, 2}, {}, {2}, {0, 2}, {1, 2}} | {{0, 1, 2}, {}, {2}, {0, 2}, {1, 2}}</td></tr>
<tr><td>44</td><td>{{0, 1, 2}, {}, {0}, {2}, {0, 2}, {1, 2}} | {{0, 1, 2}, {}, {0}, {2}, {0, 2}, {1, 2}}</td></tr>
<tr><td>45</td><td>{{0, 1, 2}, {}, {1}, {2}, {0, 2}, {1, 2}} | {{0, 1, 2}, {}, {1}, {2}, {0, 2}, {1, 2}}</td></tr>
<tr><td>46</td><td>{{0, 1, 2}, {}, {1}, {0, 2}} | {{0, 1, 2}, {}, {1}, {0, 2}}</td></tr>
</table>

| | |
|---|---|
| 47 | {{0, 1, 2}, {}, {0}, {1}, {0, 2}} \| {{0, 1, 2}, {}, {0}, {1}, {0, 2}} |
| 48 | {{0, 1, 2}, {0}, {0, 1}, {0, 2}} \| {{1, 2}, {}, {1}, {2}} |
| 49 | {{0, 1, 2}, {}, {0}, {0, 1}, {0, 2}} \| {{0, 1, 2}, {}, {0}, {0, 1}, {0, 2}} |
| 50 | {{0, 1, 2}, {}, {0}, {1}} \| {{0, 1, 2}, {}, {0}, {1}} |
| 51 | {{0, 1, 2}, {0, 1}} \| {{2}, {}} |
| 52 | {{0, 1, 2}, {}, {2}, {0, 2}} \| {{0, 1, 2}, {}, {2}, {0, 2}} |
| 53 | {{0, 1, 2}, {}, {0}, {2}, {0, 2}} \| {{0, 1, 2}, {}, {0}, {2}, {0, 2}} |
| 54 | {{0, 1, 2}, {}, {0}, {1}, {0, 1}} \| {{0, 1, 2}, {}, {0}, {1}, {0, 1}} |
| 55 | {{0, 1, 2}, {2}} \| {{0, 1}, {}} |
| 56 | {{0, 1, 2}, {}, {2}} \| {{0, 1, 2}, {}, {2}} |
| 57 | {{0, 1, 2}, {}, {0}, {2}} \| {{0, 1, 2}, {}, {0}, {2}} |
| 58 | {{0, 1, 2}, {}, {1}, {2}} \| {{0, 1, 2}, {}, {1}, {2}} |
| 59 | {{0, 1, 2}, {}, {0, 1}} \| {{0, 1, 2}, {}, {0, 1}} |
| 60 | {{0, 1, 2}, {0}, {0, 1}} \| {{1, 2}, {}, {1}} |
| 61 | {{0, 1, 2}, {}, {1}, {2}, {0, 2}} \| {{0, 1, 2}, {}, {1}, {2}, {0, 2}} |

**Table 4. Ternary non-binary convexities**

| No | G |
|---|---|
| 1 | {{0, 1, 2, 3}, {}, {0}, {1}, {2}, {1, 2}, {3}, {1, 3}, {2, 3}} |
| 2 | {{0, 1, 2, 3}, {}, {1}, {0, 1}, {2}, {1, 2}, {3}, {1, 3}, {2, 3}} |
| 3 | {{0, 1, 2, 3}, {}, {0}, {1}, {0, 1}, {2}, {1, 2}, {3}, {1, 3}, {2, 3}} |
| 4 | {{0, 1, 2, 3}, {}, {1}, {2}, {0, 2}, {1, 2}, {3}, {1, 3}, {2, 3}} |
| 5 | {{0, 1, 2, 3}, {}, {0}, {1}, {2}, {0, 2}, {1, 2}, {3}, {1, 3}, {2, 3}} |
| 6 | {{0, 1, 2, 3}, {}, {0}, {1}, {0, 1}, {2}, {0, 2}, {1, 2}, {3}, {1, 3}, {2, 3}} |
| 7 | {{0, 1, 2, 3}, {}, {0}, {1}, {2}, {1, 2}, {3}, {1, 3}, {2, 3}, {0, 2, 3}} |
| 8 | {{0, 1, 2, 3}, {}, {0}, {1}, {0, 1}, {2}, {1, 2}, {3}, {1, 3}, {2, 3}, {0, 2, 3}} |
| 9 | {{0, 1, 2, 3}, {}, {1}, {2}, {0, 2}, {1, 2}, {3}, {1, 3}, {2, 3}, {0, 2, 3}} |
| 10 | {{0, 1, 2, 3}, {}, {0}, {1}, {2}, {0, 2}, {1, 2}, {3}, {1, 3}, {2, 3}, {0, 2, 3}} |
| 11 | {{0, 1, 2, 3}, {}, {0}, {1}, {0, 1}, {2}, {0, 2}, {1, 2}, {3}, {1, 3}, {2, 3}, {1, 2, 3}} |
| 12 | {{0, 1, 2, 3}, {}, {0}, {1}, {0, 1}, {2}, {0, 2}, {1, 2}, {3}, {1, 3}, {2, 3}, {0, 2, 3}} |
| 13 | {{0, 1, 2, 3}, {}, {0}, {1}, {0, 1}, {2}, {0, 2}, {1, 2}, {3}, {1, 3}, {2, 3}, {0, 2, 3}, {1, 2, 3}} |
| 14 | {{0, 1, 2, 3}, {}, {1}, {2}, {1, 2}, {0, 1, 2}, {3}, {1, 3}, {2, 3}} |
| 15 | {{0, 1, 2, 3}, {}, {0}, {1}, {2}, {1, 2}, {0, 1, 2}, {3}, {1, 3}, {2, 3}} |
| 16 | {{0, 1, 2, 3}, {}, {1}, {0, 1}, {2}, {1, 2}, {0, 1, 2}, {3}, {1, 3}, {2, 3}} |
| 17 | {{0, 1, 2, 3}, {}, {0}, {1}, {0, 1}, {2}, {1, 2}, {0, 1, 2}, {3}, {1, 3}, {2, 3}} |
| 18 | {{0, 1, 2, 3}, {}, {1}, {2}, {0, 2}, {1, 2}, {0, 1, 2}, {3}, {1, 3}, {2, 3}, {0, 2, 3}} |
| 19 | {{0, 1, 2, 3}, {}, {0}, {1}, {2}, {0, 2}, {1, 2}, {0, 1, 2}, {3}, {1, 3}, {2, 3}, {0, 2, 3}} |
| 20 | {{0, 1, 2, 3}, {}, {0}, {1}, {0, 1}, {2}, {0, 2}, {1, 2}, {0, 1, 2}, {3}, {1, 3}, {2, 3}, {0, 2, 3}} |
| 21 | {{0, 1, 2, 3}, {}, {0}, {1}, {0, 1}, {3}, {0, 3}, {1, 3}, {2, 3}, {0, 2, 3}, {1, 2, 3}} |
| 22 | {{0, 1, 2, 3}, {}, {0}, {1}, {0, 1}, {2}, {3}, {0, 3}, {1, 3}, {2, 3}, {0, 2, 3}, {1, 2, 3}} |
| 23 | {{0, 1, 2, 3}, {}, {1}, {2}, {0, 2}, {1, 2}, {0, 1, 2}, {3}, {1, 3}, {2, 3}} |
| 24 | {{0, 1, 2, 3}, {}, {0}, {1}, {2}, {0, 2}, {1, 2}, {0, 1, 2}, {3}, {1, 3}, {2, 3}} |
| 25 | {{0, 1, 2, 3}, {}, {0}, {1}, {0, 1}, {2}, {0, 2}, {1, 2}, {0, 1, 2}, {3}, {1, 3}, {2, 3}} |
| 26 | {{0, 1, 2, 3}, {}, {0}, {1}, {0, 1}, {2}, {0, 2}, {3}, {0, 3}, {1, 3}, {2, 3}, {0, 2, 3}, {1, 2, 3}} |
| 27 | {{0, 1, 2, 3}, {}, {0}, {1}, {0, 1}, {2}, {1, 2}, {3}, {0, 3}, {1, 3}, {2, 3}, {0, 2, 3}, {1, 2, 3}} |
| 28 | {{0, 1, 2, 3}, {}, {0}, {1}, {0, 1}, {3}, {0, 3}, {1, 3}, {2, 3}, {0, 2, 3}} |
| 29 | {{0, 1, 2, 3}, {}, {0}, {1}, {0, 1}, {2}, {3}, {0, 3}, {1, 3}, {2, 3}, {0, 2, 3}} |
| 30 | {{0, 1, 2, 3}, {}, {0}, {1}, {0, 1}, {2}, {0, 2}, {1, 2}, {1, 3}, {0, 1, 3}, {1, 2, 3}} |
| 31 | {{0, 1, 2, 3}, {}, {0}, {1}, {0, 1}, {2}, {0, 2}, {3}, {0, 3}, {1, 3}, {2, 3}, {0, 2, 3}} |
| 32 | {{0, 1, 2, 3}, {}, {1}, {2}, {1, 2}, {3}, {0, 3}, {1, 3}, {2, 3}, {0, 2, 3}} |
| 33 | {{0, 1, 2, 3}, {}, {0}, {1}, {2}, {1, 2}, {3}, {0, 3}, {1, 3}, {2, 3}, {0, 2, 3}} |
| 34 | {{0, 1, 2, 3}, {}, {0}, {1}, {0, 1}, {2}, {1, 2}, {3}, {0, 3}, {1, 3}, {2, 3}, {0, 2, 3}} |
| 35 | {{0, 1, 2, 3}, {}, {0}, {1}, {0, 1}, {2}, {0, 2}, {1, 2}, {3}, {0, 3}, {0, 1, 3}, {0, 2, 3}} |
| 36 | {{0, 1, 2, 3}, {}, {0}, {1}, {0, 1}, {2}, {0, 2}, {1, 2}, {3}, {0, 3}, {1, 3}, {2, 3}, {0, 2, 3}, {1, 2, 3}} |
| 37 | {{0, 1, 2, 3}, {}, {0}, {1}, {0, 1}, {3}, {0, 3}, {1, 3}, {2, 3}, {1, 2, 3}} |
| 38 | {{0, 1, 2, 3}, {}, {0}, {1}, {0, 1}, {2}, {3}, {0, 3}, {1, 3}, {2, 3}, {1, 2, 3}} |
| 39 | {{0, 1, 2, 3}, {}, {0}, {1}, {0, 1}, {2}, {0, 2}, {1, 2}, {0, 1, 2}, {3}, {0, 3}, {1, 3}, {2, 3}, {0, 2, 3}, {1, 2, 3}} |
| 40 | {{0, 1, 2, 3}, {}, {0}, {1}, {0, 1}, {2}, {0, 2}, {1, 2}, {1, 2, 3}} |
| 41 | {{0, 1, 2, 3}, {}, {0}, {1}, {0, 1}, {3}, {0, 3}, {1, 3}, {2, 3}} |
| 42 | {{0, 1, 2, 3}, {}, {0}, {1}, {0, 1}, {2}, {3}, {0, 3}, {1, 3}, {2, 3}} |
| 43 | {{0, 1, 2, 3}, {}, {0}, {2}, {0, 2}, {3}, {0, 3}, {1, 3}, {2, 3}} |
| 44 | {{0, 1, 2, 3}, {}, {0}, {1}, {2}, {0, 2}, {3}, {0, 3}, {1, 3}, {2, 3}} |
| 45 | {{0, 1, 2, 3}, {}, {0}, {1}, {0, 1}, {2}, {0, 2}, {3}, {0, 3}, {1, 3}, {2, 3}} |
| 46 | {{0, 1, 2, 3}, {}, {1}, {2}, {1, 2}, {3}, {0, 3}, {1, 3}, {2, 3}} |
| 47 | {{0, 1, 2, 3}, {}, {0}, {1}, {2}, {1, 2}, {3}, {0, 3}, {1, 3}, {2, 3}} |
| 48 | {{0, 1, 2, 3}, {}, {0}, {1}, {0, 1}, {2}, {1, 2}, {3}, {0, 3}, {1, 3}, {2, 3}} |
| 49 | {{0, 1, 2, 3}, {}, {0}, {2}, {0, 2}, {3}, {0, 3}, {1, 3}, {2, 3}, {1, 2, 3}} |
| 50 | {{0, 1, 2, 3}, {}, {0}, {1}, {2}, {0, 2}, {3}, {0, 3}, {1, 3}, {2, 3}, {1, 2, 3}} |
| 51 | {{0, 1, 2, 3}, {}, {0}, {1}, {0, 1}, {2}, {0, 2}, {1, 2}, {3}, {0, 3}, {0, 1, 3}} |
| 52 | {{0, 1, 2, 3}, {}, {0}, {1}, {2}, {0, 2}, {1, 2}, {3}, {0, 3}, {1, 3}, {2, 3}, {0, 2, 3}} |
| 53 | {{0, 1, 2, 3}, {}, {0}, {1}, {0, 1}, {2}, {0, 2}, {1, 2}, {3}, {0, 3}, {1, 3}, {2, 3}, {0, 2, 3}} |
| 54 | {{0, 1, 2, 3}, {}, {0}, {1}, {0, 1}, {2}, {0, 2}, {1, 2}, {2, 3}, {0, 2, 3}, {1, 2, 3}} |
| 55 | {{0, 1, 2, 3}, {}, {0}, {1}, {0, 1}, {2}, {0, 2}, {3}, {0, 3}, {1, 3}, {2, 3}, {1, 2, 3}} |
| 56 | {{0, 1, 2, 3}, {}, {0}, {1}, {2}, {0, 2}, {1, 2}, {3}, {0, 3}, {1, 3}, {2, 3}} |
| 57 | {{0, 1, 2, 3}, {}, {0}, {1}, {0, 1}, {2}, {0, 2}, {0, 1, 2}, {3}, {0, 3}, {1, 3}, {2, 3}, {0, 2, 3}} |
| 58 | {{0, 1, 2, 3}, {}, {0}, {1}, {0, 1}, {2}, {0, 2}, {1, 2}, {3}, {0, 3}, {1, 3}, {2, 3}} |
| 59 | {{0, 1, 2, 3}, {}, {0}, {1}, {0, 1}, {2}, {0, 2}, {1, 2}, {0, 2, 3}} |
| 60 | {{0, 1, 2, 3}, {}, {0}, {1}, {0, 1}, {2}, {0, 1, 2}, {3}, {0, 3}, {1, 3}, {2, 3}} |
| 61 | {{0, 1, 2, 3}, {}, {0}, {1}, {0, 1}, {2}, {1, 2}, {3}, {0, 3}, {1, 3}, {2, 3}, {1, 2, 3}} |
| 62 | {{0, 1, 2, 3}, {}, {0}, {1}, {2}, {0, 2}, {1, 2}, {3}, {0, 3}, {1, 3}, {2, 3}, {1, 2, 3}} |
| 63 | {{0, 1, 2, 3}, {}, {0}, {1}, {0, 1}, {2}, {0, 2}, {1, 2}, {3}, {0, 3}, {1, 3}, {2, 3}, {1, 2, 3}} |
| 64 | {{0, 1, 2, 3}, {}, {0}, {1}, {2}, {0, 2}, {1, 2}, {0, 1, 2}, {3}, {0, 3}, {1, 3}, {2, 3}, {0, 2, 3}} |
| 65 | {{0, 1, 2, 3}, {}, {0}, {1}, {0, 1}, {2}, {0, 2}, {1, 2}, {0, 1, 2}, {3}, {0, 3}, {1, 3}, {2, 3}, {0, 2, 3}} |
| 66 | {{0, 1, 2, 3}, {}, {0}, {1}, {0, 1}, {2}, {1, 2}, {0, 1, 2}, {3}, {0, 3}, {1, 3}, {2, 3}, {1, 2, 3}} |
| 67 | {{0, 1, 2, 3}, {}, {0}, {1}, {2}, {0, 2}, {1, 2}, {0, 1, 2}, {3}, {0, 3}, {1, 3}, {2, 3}, {1, 2, 3}} |
| 68 | {{0, 1, 2, 3}, {}, {0}, {1}, {0, 1}, {2}, {0, 2}, {1, 2}, {0, 1, 2}, {3}, {0, 3}, {1, 3}, {2, 3}, {1, 2, 3}} |
| 69 | {{0, 1, 2, 3}, {}, {0}, {1}, {2}, {0, 2}, {0, 1, 2}, {3}, {0, 3}, {1, 3}, {2, 3}} |
| 70 | {{0, 1, 2, 3}, {}, {0}, {1}, {0, 1}, {2}, {0, 2}, {0, 1, 2}, {3}, {0, 3}, {1, 3}, {2, 3}} |
| 71 | {{0, 1, 2, 3}, {}, {0}, {1}, {2}, {1, 2}, {0, 1, 2}, {3}, {0, 3}, {1, 3}, {2, 3}} |
| 72 | {{0, 1, 2, 3}, {}, {0}, {1}, {0, 1}, {2}, {1, 2}, {0, 1, 2}, {3}, {0, 3}, {1, 3}, {2, 3}} |
| 73 | {{0, 1, 2, 3}, {}, {0}, {1}, {0, 1}, {2}, {0, 2}, {1, 2}, {3}, {0, 3}, {1, 3}, {0, 1, 3}, {0, 2, 3}} |
| 74 | {{0, 1, 2, 3}, {}, {0}, {1}, {0, 1}, {2}, {0, 2}, {1, 2}, {1, 3}, {0, 1, 3}} |

| 75 | {{0, 1, 2, 3}, {}, {0}, {1}, {0, 1}, {2}, {0, 2}, {1, 2}, {3}, {0, 3}, {1, 3}, {0, 1, 3}, {2, 3}, {0, 2, 3}, {1, 2, 3}} |
|---|---|
| 76 | {{0, 1, 2, 3}, {}, {0}, {1}, {0, 1}, {2}, {0, 2}, {1, 2}} |
| 77 | {{0, 1, 2, 3}, {}, {0}, {1}, {0, 1}, {2}, {0, 2}, {1, 2}, {3}, {1, 3}, {0, 1, 3}, {1, 2, 3}} |
| 78 | {{0, 1, 2, 3}, {}, {0}, {1}, {2}, {0, 2}, {1, 2}, {0, 1, 2}, {3}, {0, 3}, {1, 3}, {2, 3}} |
| 79 | {{0, 1, 2, 3}, {}, {0}, {1}, {0, 1}, {2}, {0, 2}, {1, 2}, {0, 1, 2}, {3}, {0, 3}, {1, 3}, {2, 3}} |
| 80 | {{0, 1, 2, 3}, {}, {0}, {1}, {0, 1}, {2}, {0, 2}, {1, 2}, {3}, {1, 3}, {0, 1, 3}} |
| 81 | {{0, 1, 2, 3}, {}, {0}, {1}, {0, 1}, {2}, {0, 2}, {1, 2}, {3}, {2, 3}, {0, 2, 3}, {1, 2, 3}} |
| 82 | {{0, 1, 2, 3}, {}, {0}, {1}, {0, 1}, {2}, {0, 2}, {1, 2}, {3}, {0, 3}, {0, 1, 3}, {2, 3}, {0, 2, 3}} |
| 83 | {{0, 1, 2, 3}, {}, {0}, {1}, {0, 1}, {2}, {0, 2}, {1, 2}, {3}, {0, 3}, {1, 3}, {0, 1, 3}} |
| 84 | {{0, 1, 2, 3}, {}, {0}, {1}, {0, 1}, {2}, {0, 2}, {1, 2}, {3}, {1, 2, 3}} |
| 85 | {{0, 1, 2, 3}, {}, {0}, {1}, {0, 1}, {2}, {0, 2}, {1, 2}, {2, 3}, {0, 2, 3}} |
| 86 | {{0, 1, 2, 3}, {}, {0}, {1}, {0, 1}, {2}, {0, 2}, {1, 2}, {3}, {0, 3}, {1, 3}, {0, 1, 3}, {1, 2, 3}} |
| 87 | {{0, 1, 2, 3}, {}, {0}, {1}, {0, 1}, {2}, {0, 2}, {1, 2}, {3}, {0, 2, 3}} |
| 88 | {{0, 1, 2, 3}, {}, {0}, {1}, {0, 1}, {2}, {0, 2}, {1, 2}, {3}, {0, 3}, {2, 3}, {0, 2, 3}, {1, 2, 3}} |
| 89 | {{0, 1, 2, 3}, {}, {1}, {2}, {1, 2}, {3}, {0, 3}, {1, 3}, {0, 1, 3}, {2, 3}, {0, 2, 3}} |
| 90 | {{0, 1, 2, 3}, {}, {0}, {1}, {2}, {1, 2}, {3}, {0, 3}, {1, 3}, {0, 1, 3}, {2, 3}, {0, 2, 3}} |
| 91 | {{0, 1, 2, 3}, {}, {0}, {1}, {0, 1}, {2}, {1, 2}, {3}, {0, 3}, {1, 3}, {0, 1, 3}, {2, 3}, {0, 2, 3}} |
| 92 | {{0, 1, 2, 3}, {}, {0}, {1}, {2}, {0, 2}, {1, 2}, {3}, {0, 3}, {1, 3}, {0, 1, 3}, {2, 3}, {0, 2, 3}} |
| 93 | {{0, 1, 2, 3}, {}, {0}, {1}, {0, 1}, {2}, {0, 2}, {1, 2}, {3}, {1, 3}, {0, 1, 3}, {2, 3}, {1, 2, 3}} |
| 94 | {{0, 1, 2, 3}, {}, {0}, {1}, {0, 1}, {2}, {0, 2}, {1, 2}, {3}, {0, 3}, {1, 3}, {0, 1, 3}, {2, 3}, {0, 2, 3}} |
| 95 | {{0, 1, 2, 3}, {}, {0}, {1}, {0, 1}, {2}, {0, 2}, {1, 2}, {3}, {2, 3}, {0, 2, 3}} |
| 96 | {{0, 1, 2, 3}, {}, {0}, {1}, {0, 1}, {2}, {0, 2}, {1, 2}, {3}} |
| 97 | {{0, 1, 2, 3}, {}, {0}, {1}, {0, 1}, {2}, {0, 2}, {1, 2}, {3}, {0, 1, 3}, {2, 3}} |
| 98 | {{0, 1, 2, 3}, {}, {0}, {2}, {0, 2}, {3}, {0, 3}, {1, 3}, {0, 1, 3}, {2, 3}, {1, 2, 3}} |
| 99 | {{0, 1, 2, 3}, {}, {0}, {1}, {2}, {0, 2}, {3}, {0, 3}, {1, 3}, {0, 1, 3}, {2, 3}, {1, 2, 3}} |
| 100 | {{0, 1, 2, 3}, {}, {0}, {1}, {0, 1}, {2}, {0, 2}, {3}, {0, 3}, {1, 3}, {0, 1, 3}, {2, 3}, {1, 2, 3}} |
| 101 | {{0, 1, 2, 3}, {}, {0}, {1}, {0, 1}, {2}, {0, 2}, {1, 2}, {0, 1, 2}, {3}, {0, 3}, {1, 3}, {0, 1, 3}, {2, 3}, {0, 2, 3}} |
| 102 | {{0, 1, 2, 3}, {}, {0}, {1}, {0, 1}, {2}, {0, 2}, {1, 2}, {3}, {0, 3}, {1, 2, 3}} |
| 103 | {{0, 1, 2, 3}, {}, {0}, {1}, {0, 1}, {2}, {0, 2}, {1, 2}, {2, 3}, {1, 2, 3}} |
| 104 | {{0, 1, 2, 3}, {}, {0}, {1}, {0, 1}, {2}, {0, 2}, {1, 2}, {0, 3}, {0, 2, 3}} |
| 105 | {{0, 1, 2, 3}, {}, {0}, {2}, {0, 2}, {3}, {0, 3}, {0, 1, 3}, {2, 3}} |
| 106 | {{0, 1, 2, 3}, {}, {0}, {1}, {2}, {0, 2}, {3}, {0, 3}, {0, 1, 3}, {2, 3}} |
| 107 | {{0, 1, 2, 3}, {}, {0}, {0, 1}, {2}, {0, 2}, {3}, {0, 3}, {0, 1, 3}, {2, 3}} |
| 108 | {{0, 1, 2, 3}, {}, {0}, {1}, {0, 1}, {2}, {0, 2}, {3}, {0, 3}, {0, 1, 3}, {2, 3}} |
| 109 | {{0, 1, 2, 3}, {}, {0}, {1}, {2}, {0, 2}, {1, 2}, {3}, {0, 3}, {1, 3}, {0, 1, 3}, {2, 3}, {1, 2, 3}} |
| 110 | {{0, 1, 2, 3}, {}, {0}, {1}, {0, 1}, {2}, {0, 2}, {1, 2}, {3}, {0, 3}, {1, 3}, {0, 1, 3}, {2, 3}, {1, 2, 3}} |
| 111 | {{0, 1, 2, 3}, {}, {0}, {1}, {0, 1}, {2}, {0, 2}, {1, 2}, {0, 1, 2}, {3}, {0, 3}, {1, 3}, {0, 1, 3}, {2, 3}, {1, 2, 3}} |
| 112 | {{0, 1, 2, 3}, {}, {0}, {1}, {0, 1}, {2}, {0, 2}, {1, 2}, {3}, {0, 3}, {2, 3}, {0, 2, 3}} |
| 113 | {{0, 1, 2, 3}, {}, {1}, {2}, {1, 2}, {3}, {1, 3}, {2, 3}, {0, 2, 3}} |
| 114 | {{0, 1, 2, 3}, {}, {0}, {1}, {0, 1}, {2}, {0, 2}, {1, 2}, {1, 3}, {1, 2, 3}} |
| 115 | {{0, 1, 2, 3}, {}, {0}, {1}, {2}, {0, 2}, {1, 2}, {3}, {0, 3}, {0, 1, 3}, {2, 3}} |
| 116 | {{0, 1, 2, 3}, {}, {0}, {1}, {0, 1}, {2}, {0, 2}, {1, 2}, {3}, {0, 3}, {0, 1, 3}, {2, 3}} |
| 117 | {{0, 1, 2, 3}, {}, {0}, {1}, {0, 1}, {2}, {0, 2}, {1, 2}, {3}, {1, 3}, {1, 2, 3}} |
| 118 | {{0, 1, 2, 3}, {}, {0}, {1}, {0, 1}, {2}, {0, 2}, {1, 2}, {3}, {2, 3}, {1, 2, 3}} |
| 119 | {{0, 1, 2, 3}, {}, {0}, {0, 1}, {2}, {0, 2}, {0, 1, 2}, {3}, {0, 3}, {0, 1, 3}, {2, 3}} |
| 120 | {{0, 1, 2, 3}, {}, {0}, {1}, {0, 1}, {2}, {0, 2}, {0, 1, 2}, {3}, {0, 3}, {0, 1, 3}, {2, 3}} |
| 121 | {{0, 1, 2, 3}, {}, {0}, {1}, {0, 1}, {2}, {0, 2}, {1, 2}, {3}, {0, 3}, {0, 2, 3}} |
| 122 | {{0, 1, 2, 3}, {}, {0}, {1}, {0, 1}, {2}, {0, 2}, {1, 2}, {2, 3}} |
| 123 | {{0, 1, 2, 3}, {}, {0}, {1}, {0, 1}, {2}, {0, 2}, {1, 2}, {0, 3}} |
| 124 | {{0, 1, 2, 3}, {}, {0}, {1}, {0, 1}, {3}, {0, 3}, {1, 3}, {1, 2, 3}} |
| 125 | {{0, 1, 2, 3}, {}, {0}, {1}, {0, 1}, {2}, {3}, {0, 3}, {1, 3}, {1, 2, 3}} |
| 126 | {{0, 1, 2, 3}, {}, {0}, {1}, {0, 1}, {2}, {0, 2}, {3}, {0, 3}, {1, 3}, {1, 2, 3}} |
| 127 | {{0, 1, 2, 3}, {}, {0}, {1}, {0, 1}, {2}, {0, 2}, {1, 2}, {0, 1, 2}, {3}, {0, 3}, {0, 1, 3}, {2, 3}} |
| 128 | {{0, 1, 2, 3}, {}, {0}, {1}, {0, 1}, {1, 2}, {3}, {0, 3}, {1, 3}, {1, 2, 3}} |
| 129 | {{0, 1, 2, 3}, {}, {0}, {1}, {0, 1}, {2}, {0, 2}, {1, 2}, {3}, {0, 3}} |
| 130 | {{0, 1, 2, 3}, {}, {0}, {1}, {0, 1}, {2}, {0, 2}, {1, 2}, {3}, {1, 3}, {0, 2, 3}} |
| 131 | {{0, 1, 2, 3}, {}, {0}, {2}, {0, 2}, {3}, {0, 3}, {2, 3}, {1, 2, 3}} |
| 132 | {{0, 1, 2, 3}, {}, {0}, {1}, {0, 1}, {2}, {1, 2}, {3}, {0, 3}, {1, 3}, {1, 2, 3}} |
| 133 | {{0, 1, 2, 3}, {}, {0}, {1}, {2}, {0, 2}, {3}, {0, 3}, {2, 3}, {1, 2, 3}} |
| 134 | {{0, 1, 2, 3}, {}, {0}, {1}, {0, 1}, {2}, {0, 2}, {3}, {0, 3}, {2, 3}, {1, 2, 3}} |
| 135 | {{0, 1, 2, 3}, {}, {0}, {1}, {0, 1}, {2}, {0, 2}, {1, 2}, {3}, {0, 3}, {1, 3}, {1, 2, 3}} |
| 136 | {{0, 1, 2, 3}, {}, {1}, {2}, {1, 2}, {3}, {1, 3}, {0, 1, 3}, {2, 3}} |
| 137 | {{0, 1, 2, 3}, {}, {0}, {1}, {2}, {1, 2}, {3}, {1, 3}, {0, 1, 3}, {2, 3}} |
| 138 | {{0, 1, 2, 3}, {}, {1}, {0, 1}, {2}, {1, 2}, {3}, {1, 3}, {0, 1, 3}, {2, 3}} |
| 139 | {{0, 1, 2, 3}, {}, {0}, {1}, {0, 1}, {2}, {1, 2}, {3}, {1, 3}, {0, 1, 3}, {2, 3}} |
| 140 | {{0, 1, 2, 3}, {}, {0}, {1}, {0, 1}, {2}, {0, 2}, {1, 2}, {1, 3}} |
| 141 | {{0, 1, 2, 3}, {}, {0}, {2}, {0, 2}, {1, 2}, {3}, {0, 3}, {2, 3}, {1, 2, 3}} |
| 142 | {{0, 1, 2, 3}, {}, {0}, {1}, {2}, {0, 2}, {1, 2}, {3}, {0, 3}, {2, 3}, {1, 2, 3}} |
| 143 | {{0, 1, 2, 3}, {}, {0}, {1}, {0, 1}, {2}, {0, 2}, {1, 2}, {3}, {0, 3}, {2, 3}, {1, 2, 3}} |
| 144 | {{0, 1, 2, 3}, {}, {0}, {1}, {0, 1}, {3}, {0, 3}, {1, 3}, {0, 2, 3}} |
| 145 | {{0, 1, 2, 3}, {}, {0}, {1}, {0, 1}, {2}, {3}, {0, 3}, {1, 3}, {0, 2, 3}} |
| 146 | {{0, 1, 2, 3}, {}, {0}, {1}, {0, 1}, {0, 2}, {3}, {0, 3}, {1, 3}, {0, 2, 3}} |
| 147 | {{0, 1, 2, 3}, {}, {0}, {1}, {0, 1}, {2}, {0, 2}, {3}, {0, 3}, {1, 3}, {0, 2, 3}} |
| 148 | {{0, 1, 2, 3}, {}, {0}, {1}, {0, 1}, {1, 2}, {0, 1, 2}, {3}, {0, 3}, {1, 3}, {1, 2, 3}} |
| 149 | {{0, 1, 2, 3}, {}, {0}, {1}, {0, 1}, {2}, {1, 2}, {0, 1, 2}, {3}, {0, 3}, {1, 3}, {1, 2, 3}} |
| 150 | {{0, 1, 2, 3}, {}, {0}, {1}, {0, 1}, {2}, {0, 2}, {1, 2}, {0, 1, 2}, {3}, {0, 3}, {1, 3}, {1, 2, 3}} |
| 151 | {{0, 1, 2, 3}, {}, {0}, {1}, {2}, {0, 2}, {1, 2}, {3}, {1, 3}, {0, 1, 3}, {2, 3}} |
| 152 | {{0, 1, 2, 3}, {}, {0}, {1}, {0, 1}, {2}, {0, 2}, {1, 2}, {3}, {1, 3}, {0, 1, 3}, {2, 3}} |
| 153 | {{0, 1, 2, 3}, {}, {0}, {2}, {0, 2}, {1, 2}, {0, 1, 2}, {3}, {0, 3}, {2, 3}, {1, 2, 3}} |
| 154 | {{0, 1, 2, 3}, {}, {0}, {1}, {2}, {0, 2}, {1, 2}, {0, 1, 2}, {3}, {0, 3}, {2, 3}, {1, 2, 3}} |
| 155 | {{0, 1, 2, 3}, {}, {1}, {0, 1}, {2}, {1, 2}, {0, 1, 2}, {3}, {1, 3}, {0, 1, 3}, {2, 3}} |
| 156 | {{0, 1, 2, 3}, {}, {0}, {1}, {0, 1}, {2}, {1, 2}, {0, 1, 2}, {3}, {1, 3}, {0, 1, 3}, {2, 3}} |
| 157 | {{0, 1, 2, 3}, {}, {0}, {1}, {0, 1}, {2}, {0, 2}, {1, 2}, {0, 1, 2}, {3}, {0, 3}, {2, 3}, {1, 2, 3}} |
| 158 | {{0, 1, 2, 3}, {}, {0}, {1}, {0, 1}, {2}, {1, 2}, {3}, {0, 3}, {1, 3}, {0, 2, 3}} |
| 159 | {{0, 1, 2, 3}, {}, {0}, {1}, {0, 1}, {2}, {0, 2}, {1, 2}, {0, 1, 2}, {3}, {1, 3}, {0, 1, 3}, {2, 3}} |
| 160 | {{0, 1, 2, 3}, {}, {0}, {1}, {0, 1}, {2}, {0, 2}, {1, 2}, {3}, {0, 3}, {1, 3}, {0, 2, 3}} |
| 161 | {{0, 1, 2, 3}, {}, {0}, {1}, {0, 1}, {0, 2}, {0, 1, 2}, {3}, {0, 3}, {1, 3}, {0, 2, 3}} |
| 162 | {{0, 1, 2, 3}, {}, {0}, {1}, {0, 1}, {2}, {0, 2}, {1, 2}, {3}, {2, 3}} |
| 163 | {{0, 1, 2, 3}, {}, {0}, {1}, {0, 1}, {2}, {0, 2}, {1, 2}, {3}, {1, 3}} |
| 164 | {{0, 1, 2, 3}, {}, {0}, {1}, {0, 1}, {2}, {0, 2}, {0, 1, 2}, {3}, {0, 3}, {1, 3}, {0, 2, 3}} |
| 165 | {{0, 1, 2, 3}, {}, {0}, {2}, {0, 2}, {3}, {0, 3}, {1, 3}, {0, 1, 3}, {2, 3}} |

| 166 | {{0, 1, 2, 3}, {}, {0}, {1}, {2}, {0, 2}, {3}, {0, 3}, {1, 3}, {0, 1, 3}, {2, 3}} |
|---|---|
| 167 | {{0, 1, 2, 3}, {}, {0}, {1}, {0, 1}, {2}, {0, 2}, {3}, {0, 3}, {1, 3}, {0, 1, 3}, {2, 3}} |
| 168 | {{0, 1, 2, 3}, {}, {1}, {2}, {1, 2}, {3}, {0, 3}, {1, 3}, {0, 1, 3}, {2, 3}} |
| 169 | {{0, 1, 2, 3}, {}, {0}, {1}, {2}, {1, 2}, {3}, {0, 3}, {1, 3}, {0, 1, 3}, {2, 3}} |
| 170 | {{0, 1, 2, 3}, {}, {0}, {1}, {0, 1}, {2}, {0, 2}, {1, 2}, {0, 1, 2}, {3}, {0, 3}, {1, 3}, {0, 2, 3}} |
| 171 | {{0, 1, 2, 3}, {}, {0}, {1}, {0, 1}, {3}, {0, 3}, {1, 3}} |
| 172 | {{0, 1, 2, 3}, {}, {0}, {1}, {0, 1}, {2}, {3}, {0, 3}, {1, 3}} |
| 173 | {{0, 1, 2, 3}, {}, {0}, {1}, {0, 1}, {0, 2}, {3}, {0, 3}, {1, 3}} |
| 174 | {{0, 1, 2, 3}, {}, {0}, {1}, {0, 1}, {2}, {0, 2}, {3}, {0, 3}, {1, 3}} |
| 175 | {{0, 1, 2, 3}, {}, {0}, {1}, {0, 1}, {1, 2}, {3}, {0, 3}, {1, 3}} |
| 176 | {{0, 1, 2, 3}, {}, {0}, {1}, {0, 1}, {2}, {1, 2}, {3}, {0, 3}, {1, 3}} |
| 177 | {{0, 1, 2, 3}, {}, {0}, {1}, {0, 1}, {2}, {1, 2}, {3}, {0, 3}, {1, 3}, {0, 1, 3}, {2, 3}} |
| 178 | {{0, 1, 2, 3}, {}, {0}, {1}, {2}, {0, 2}, {1, 2}, {3}, {0, 3}, {1, 3}, {0, 1, 3}, {2, 3}} |
| 179 | {{0, 1, 2, 3}, {}, {0}, {1}, {0, 1}, {2}, {0, 2}, {1, 2}, {3}, {0, 3}, {1, 3}} |
| 180 | {{0, 1, 2, 3}, {}, {0}, {1}, {0, 1}, {0, 1, 2}, {3}, {0, 3}, {1, 3}} |
| 181 | {{0, 1, 2, 3}, {}, {0}, {1}, {0, 1}, {2}, {0, 2}, {1, 2}, {3}, {0, 3}, {1, 3}, {0, 1, 3}, {2, 3}} |
| 182 | {{0, 1, 2, 3}, {}, {0}, {1}, {0, 1}, {2}, {0, 1, 2}, {3}, {0, 3}, {1, 3}} |
| 183 | {{0, 1, 2, 3}, {}, {0}, {1}, {0, 1}, {0, 2}, {0, 1, 2}, {3}, {0, 3}, {1, 3}} |
| 184 | {{0, 1, 2, 3}, {}, {0}, {1}, {0, 1}, {2}, {0, 2}, {0, 1, 2}, {3}, {0, 3}, {1, 3}} |
| 185 | {{0, 1, 2, 3}, {}, {0}, {1}, {0, 1}, {2}, {0, 2}, {0, 1, 2}, {3}, {0, 3}, {1, 3}, {0, 1, 3}, {2, 3}} |
| 186 | {{0, 1, 2, 3}, {}, {0}, {1}, {0, 1}, {1, 2}, {0, 1, 2}, {3}, {0, 3}, {1, 3}} |
| 187 | {{0, 1, 2, 3}, {}, {0}, {1}, {0, 1}, {2}, {1, 2}, {0, 1, 2}, {3}, {0, 3}, {1, 3}, {0, 1, 3}, {2, 3}} |
| 188 | {{0, 1, 2, 3}, {}, {0}, {1}, {0, 1}, {2}, {0, 2}, {1, 2}, {0, 1, 2}, {3}, {0, 3}, {1, 3}, {0, 1, 3}, {2, 3}} |
| 189 | {{0, 1, 2, 3}, {}, {0}, {1}, {0, 1}, {2}, {1, 2}, {0, 1, 2}, {3}, {0, 3}, {1, 3}} |
| 190 | {{0, 1, 2, 3}, {}, {0}, {1}, {0, 1}, {2}, {0, 2}, {1, 2}, {0, 1, 2}, {3}, {0, 3}, {1, 3}} |
| 191 | {{0, 1, 2, 3}, {}, {0}, {2}, {0, 2}, {3}, {0, 3}, {2, 3}} |
| 192 | {{0, 1, 2, 3}, {}, {0}, {1}, {2}, {0, 2}, {3}, {0, 3}, {2, 3}} |
| 193 | {{0, 1, 2, 3}, {}, {0}, {0, 1}, {2}, {0, 2}, {3}, {0, 3}, {2, 3}} |
| 194 | {{0, 1, 2, 3}, {}, {0}, {1}, {0, 1}, {2}, {0, 2}, {3}, {0, 3}, {2, 3}} |
| 195 | {{0, 1, 2, 3}, {}, {0}, {1}, {0, 1}, {2}, {0, 2}, {1, 2}, {0, 3}, {0, 1, 3}, {0, 2, 3}} |
| 196 | {{0, 1, 2, 3}, {}, {0}, {2}, {0, 2}, {1, 2}, {3}, {0, 3}, {2, 3}} |
| 197 | {{0, 1, 2, 3}, {}, {0}, {1}, {2}, {0, 2}, {1, 2}, {3}, {0, 3}, {2, 3}} |
| 198 | {{0, 1, 2, 3}, {}, {0}, {1}, {0, 1}, {2}, {0, 2}, {1, 2}, {3}, {0, 3}, {2, 3}} |
| 199 | {{0, 1, 2, 3}, {}, {0}, {2}, {0, 2}, {0, 1, 2}, {3}, {0, 3}, {2, 3}} |
| 200 | {{0, 1, 2, 3}, {}, {0}, {1}, {2}, {0, 2}, {0, 1, 2}, {3}, {0, 3}, {2, 3}} |
| 201 | {{0, 1, 2, 3}, {}, {0}, {0, 1}, {2}, {0, 2}, {0, 1, 2}, {3}, {0, 3}, {2, 3}} |
| 202 | {{0, 1, 2, 3}, {}, {0}, {1}, {0, 1}, {2}, {0, 2}, {0, 1, 2}, {3}, {0, 3}, {2, 3}} |
| 203 | {{0, 1, 2, 3}, {}, {0}, {2}, {0, 2}, {1, 2}, {0, 1, 2}, {3}, {0, 3}, {2, 3}} |
| 204 | {{0, 1, 2, 3}, {}, {0}, {1}, {2}, {0, 2}, {1, 2}, {0, 1, 2}, {3}, {0, 3}, {2, 3}} |
| 205 | {{0, 1, 2, 3}, {}, {0}, {1}, {0, 1}, {2}, {0, 2}, {1, 2}, {0, 1, 3}} |
| 206 | {{0, 1, 2, 3}, {}, {0}, {1}, {0, 1}, {2}, {0, 2}, {1, 2}, {0, 1, 2}, {3}, {0, 3}, {2, 3}} |
| 207 | {{0, 1, 2, 3}, {}, {0}, {1}, {0, 1}, {2}, {0, 2}, {1, 2}, {3}, {0, 1, 3}} |
| 208 | {{0, 1, 2, 3}, {}, {1}, {2}, {1, 2}, {3}, {1, 3}, {2, 3}} |
| 209 | {{0, 1, 2, 3}, {}, {0}, {1}, {0, 1}, {2}, {0, 2}, {1, 2}, {0, 3}, {0, 1, 3}} |

# БИНОМИАЛЬНОЕ ПРЕОБРАЗОВАНИЕ ПОСЛЕДОВАТЕЛЬНОСТЕЙ, ПЕРЕЧИСЛЯЮЩИХ *N*-АРНЫЕ ВЫПУКЛОСТИ

А.М. Дуллиев [*], Д.Д. Наумихин [*]

8 июня 2026 г.

## Аннотация

Настоящее исследование посвящено перечислению $N$-арных выпуклых структур на конечных носителях. В работе получено соотношение, связывающее общее число выпуклостей с числом заземлённых выпуклостей; оно выражается в виде биномиального преобразования соответствующих последовательностей. Для $|X| \leq 5$ приведены точные значения количеств всех $N$-арных и заземлённых $N$-арных выпуклостей. Полученные целочисленные последовательности сопоставлены с данными OEIS; обнаружены как известные последовательности (например, A000798), так и новые, не зарегистрированные ранее.

[*] Кафедра прикладной математики и информатики, Казанский национально-исследовательский технический университет им. А.Н. Туполева, Российская Федерация. E-mail: {dulliev@yandex.ru, dan.naumihin@yandex.ru}.

# 1. ВВЕДЕНИЕ

Пусть дано произвольное множество $X$. *Выпуклостью* $\mathbf{G}$ в $X$ называется семейство подмножеств в $X$, замкнутое относительно произвольных пересечений своих элементов и содержащее $X$; при этом пара $(X, \mathbf{G})$ называется *выпуклым пространством или пространством выпуклости (пространством с выпуклостью)*. Каждой выпуклости взаимно однозначно соответствует *оператор выпуклой оболочки* $g: \mathcal{P}(X) \to \mathcal{P}(X)$, удовлетворяющий для любых $A \subset B \subset X$ аксиомам: $A \subset gA$, $gA \subset gB$, $ggA = gA$, — при этом соответствие устанавливается правилами: $gA = \bigcap_{A \subset B \in \mathbf{G}} B$ и $\mathbf{G} = \{A \subset X: gA = A\}$. Выпуклость $\mathbf{G}$ называется $N$-арной, если $A \in \mathbf{G} \Leftrightarrow \forall B \subset A: (|B| \leq N) \Rightarrow gB \subset A$ при $A \subset X$. Будем говорить, что выпуклость $\mathbf{G}$ удовлетворяет аксиоме заземления (такую выпуклость будем называть также выпуклостью с заземлением или заземленной “grounded convexity”), если выполняется условие $\emptyset \in \mathbf{G}$.

Систематическому изложению основ теории выпуклости посвящены две, уже ставшие классическими, монографии (В.П. Солтана [1] и M.L.J. van de Vel [2]). Стоит отметить, что определения выпуклых пространств в них различны (неэквивалентны): в [1] приводится более общее определение, а в [2] к нему добавляется условие индуктивности выпуклости. Вопросы, связанные с перечислением выпуклых структур различных типов и структур, порожденных выпуклыми, на конечных носителях, изучены весьма слабо. Результаты, полученные на сегодняшний день, в основном сводятся к перечислению (неполному) других известных структур, для которых доказана их эквивалентность некоторым выпуклым структурам и проведено перечисление. Совокупность чисел, каждое из которых равно количеству различных в том или ином смысле выпуклых структур на конечном носителе и зависит от числа его элементов, формирует естественным образом целочисленную последовательность. Иногда удается найти замкнутую формулу для произвольного члена такой последовательности. В настоящей работе получено соотношение, связывающее общее число выпуклостей с числом заземлённых выпуклостей; оно выражается в виде биномиального преобразования соответствующих последовательностей.

Целочисленные последовательности, полученные по результатам работы разработанной компьютерной программы, сравниваются с последовательностями из международной “Онлайн энциклопедии целочисленных последовательностей”. Несколько последовательностей, как и ожидалось, оказались новыми.

Приведём список символов и обозначений, используемых в работе:

- $|A|$ — мощность множества $A$;
- $\mathcal{P}(X)$ — семейство всех подмножеств множества $X$, буллеан множества $X$;
- $\mathbf{G}$ — выпуклость в $X$;
- $gA = \cap \{B \in \mathbf{G}: A \subset B\}$ — выпуклая оболочка множества $A \subset X$
- $\Gamma^N(X)$ — семейство всех $N$-арных выпуклостей на носителе $X$;
- $\Gamma_0^N(X)$ — семейство всех заземлённых $N$-арных выпуклостей на носителе $X$;
- $X_n$ — конечное множество мощности $n$.

# 2. СВОЙСТВА $\Gamma^N(X)$, $\Gamma_0^N(X)$

Всюду ниже считается, что $X$ – конечное множество.

**Лемма 1.** *Пусть* $\mathbf{G} \in \Gamma^N(X)$, $C$ *— наименьшее выпуклое множество*

$$C = g\emptyset = \bigcap_{A \in \mathbf{G}} A,$$

*а множество* $\mathbf{H}$ *построено* по *правилу*

$$\mathbf{H} = \{A \backslash C : A \in \mathbf{G}\}. \tag{1}$$

*Тогда между* $\mathbf{G}$ *и* $\mathbf{H}$ *существует взаимооднозначное соответствие.*

**Доказательство.** Из построения видно, что отображение $\varphi : X \mapsto X \backslash C$ сюръективно:

$$\forall B \in \mathbf{H} \exists A \in \mathbf{G}: B = A \backslash C.$$

Покажем инъективность $\varphi$. Предположим обратное: найдутся два различных множества $A_1, A_2$, для которых $\varphi A_1 = \varphi A_2$, т.е.

$$A_1 \neq A_2, A_1 \backslash C = A_2 \backslash C. \tag{2}$$

Из построения $C$:

$$\forall A \in \mathbf{G}: C \subset A. \tag{3}$$

Рассмотрим произвольную точку $x$ из $A_1$. Пусть $x$ находится и в $C$, тогда по (3) $x$ принадлежит и $A_2$. В противном случае, если $x \in A_1 \backslash C$, то по (2) $x$ принадлежит и $A_2 \backslash C$. Тогда $A_1 \subset A_2$. Проводя аналогичные рассуждения уже для $A_2$, получаем $A_2 \subset A_1$. И таким образом, $A_1 = A_2$, что противоречит построению (2).

Отображение $X \mapsto X \backslash C$ и инъективно, и сюръективно, а значит взаимооднозначно. ■

**Лемма 2.** *Пусть* $\mathbf{G} \in \Gamma^N(X)$, *тогда* $\mathbf{H}$, *построенное по правилу (1), является заземленной выпуклостью на* $Y = X \backslash C$.

**Доказательство.** 1. Очевидно, что согласно (1) имеем: $Y \in \mathbf{H}$.

2. Это семейство множеств замкнуто относительно пересечения

$$\forall B_1, B_2 \in \mathbf{H} \exists A_1, A_2 \in \mathbf{G}: B_1 = A_1 \backslash C, B_2 = A_2 \backslash C;$$
$$A_1 \cap A_2 \in \mathbf{G} \Rightarrow (A_1 \cap A_2) \backslash C = B_1 \cap B_2 \in \mathbf{H}.$$

3. Данное семейство множеств содержит пустое по построению:

$$g\emptyset \in G \Rightarrow g\emptyset \backslash g\emptyset = \emptyset \in \mathbf{H}.$$ ■

**Лемма 3.** *Пусть* $\mathbf{G} \in \Gamma^N(X)$. *Выпуклая оболочка* $h$ *выпуклости* $\mathbf{H}$, *построенной по правилу (1), имеет следующую связь с выпуклой оболочкой* $g$:

$$\forall A \subset Y: gA = C \cup hA \tag{4}$$

*или*

$$hA = gA \backslash C. \tag{5}$$

*К тому же, имеет место представление*

$$A \subset X \Rightarrow gA = h(A \backslash C) \cup C. \tag{6}$$

**Доказательство.** Справедливость соотношения (4) вытекает из следующей цепочки равенств:

$$gA = \bigcap_{\substack{B \in \mathbf{G} \\ A \subset B}} B = \bigcap_{\substack{Z \cup C \in \mathbf{G} \\ A \subset Z \cup C}} (Z \cup C) =$$
$$= C \cup \bigcap_{\substack{Z \cup C \in \mathbf{G} \\ A \subset Z \cup C}} Z = C \cup \bigcap_{\substack{Z \in \mathbf{H} \\ A \subset Z}} Z = C \cup hA.$$

По построению выпуклой оболочки $h$ имеет место:

$$hA \in \mathbf{H}, \mathbf{H} = \{B \backslash C : B \in \mathbf{G}\}.$$

отсюда $hA \cap C = \emptyset$, а значит, из уже доказанного выше следует соотношение (5).

Для доказательства соотношения (6) рассмотрим произвольное $A \subset X$. Заметим, что $gA = g(A \backslash C)$. Действительно, добавление к $A$ элементов из $C$ не может расширить выпуклую оболочку, так как $C = g\emptyset \subset g(A \backslash C)$ (минимальное выпуклое множество $C$ содержится в любой выпуклой оболочке). Следовательно, для любого $x \in A \cap C$ имеем $x \in C \subset g(A \backslash C)$, а значит, $gA = g(A \backslash C)$.

Теперь $A \backslash C \subset Y = X \backslash C$, потому можно применить формулу (4) к множеству $A \backslash C$:

$$g(A \backslash C) = C \cup h(A \backslash C).$$

Таким образом,

$$gA = C \cup h(A \backslash C),$$

что в точности соответствует формуле (6). ■

**Лемма 4.** *Пусть* $\mathbf{G} \in \Gamma^N(X)$. *Выпуклость* $\mathbf{H}$*, построенная по правилу (1), является* $N$*-арной.*

**Доказательство.** Рассмотрим произвольную выпуклость $\mathbf{G}$ из $\Gamma^N(X)$. Тогда по определению $N$-арной выпуклости, для неё имеем:

$$\forall A \subset X : [A \in \mathbf{G} \Leftrightarrow (\forall B \subset A : |B| \le N \Rightarrow gB \subset A)].$$

Из леммы 2 и связи выпуклых оболочек $g$ и $h$ (лемма 3), сузив исследуемое свойство до множества $Y = X \backslash C$:

$$\forall A \subset Y : [A \cup C \in \mathbf{G} \Leftrightarrow (\forall B \subset C \cup A : |B| \le N \Rightarrow hB \cup C \subset A \cup C)].$$

Применив отображение $\mathbf{G}$ на $\mathbf{H}$, а также взяв во внимание само построение выпуклой оболочки $h$, а именно тот факт, что $\mathbf{H}$ построено как подмножество $X \backslash C$, выходит, что всякая выпуклая оболочка $h$ содержится полностью в $X \backslash C$, ровно как и само множество $A$, принадлежащее $\mathbf{H}$ в левой части, то есть справедливо следующее:

$$\forall A \subset Y [A \in \mathbf{H} \Leftrightarrow (\forall B \subset A : |B| \le N \Rightarrow hB \subset A)].$$

А это соответствует определению $N$-арного множества на $Y = X \backslash C$. Из доказанной леммы 2 видно, что эта выпуклость и заземленная, так что $\mathbf{H}$ действительно заземленная выпуклость на носителе $Y = X \backslash C$. ■

Пусть $C \subset X$ — произвольное. Определим через $E_C$ семейство выпуклостей из $\Gamma^N(X)$ следующим образом:

$$\mathbf{G} \in E_C \Leftrightarrow g\emptyset = C.$$

Имеет место следующая лемма.

**Лемма 5.** *Каждой выпуклости* $\mathbf{G} \in E_C$ *можно поставить в соответствие единственную выпуклость* $\boldsymbol{H} \in \Gamma_0^N(X \backslash C)$*, задаваемую правилом (1).*

*Обратно, для любой* $\mathbf{H} \in \Gamma_0^N(X \backslash C)$ *семейство*

$$\mathbf{G} = \{C \cup D : D \in \mathbf{H}\}$$

*принадлежит* $E_C$.

*Построенное соответствие устанавливает между* $E_C$ *и* $\Gamma_0^N(X \backslash C)$ *биекцию* $\xi_C$.

**Доказательство.** (*прямое отображение*): Пусть $\mathbf{G} \in E_C$. По условию $g\emptyset = C$. Построим $\mathbf{H}$ по правилу (1). Леммы 2 и 3 гарантируют, что $\mathbf{H} \in \Gamma_0^N(X \backslash C)$. Лемма 1 гарантирует, что построенное отображение биекция.

(*обратное отображение*): Пусть $\mathbf{H} \in \Gamma_0^N(Y)$,где $Y = X \backslash C$. Определим $\mathbf{G} = \{C \cup D : D \in \mathbf{H}\}$. Проверим, что $\mathbf{G} \in E_C$:

1. По построению $C = C \cup \emptyset \in \mathbf{G}$, так как $\emptyset \in \mathbf{H}$. Более того, для любого $C \cup D \in \mathbf{G}$ имеем $C \subset C \cup D$. Следовательно, $C$ содержится во всех элементах $\mathbf{G}$, а значит, $C \subset g\emptyset$.

2. Покажем, что $g\emptyset \subset C$. Пусть $x \in g\emptyset$. Тогда $x$ принадлежит всем элементам **G**, в том числе $C$. Поэтому $x \in C$.

3. Из предыдущих пунктов $g\emptyset = C$.

4. Проверим, что **G** — $N$-арная выпуклость:

(a) **G** замкнуто относительно пересечений. Для любых $C \cup D_1, C \cup D_2 \in \mathbf{G}$:

$$(C \cup D_1) \cap (C \cup D_2) = C \cup (D_1 \cap D_2).$$

Поскольку $D_1, D_2 \in \mathbf{H}$, $D_1 \cap D_2 \in \mathbf{H}$. Поэтому $C \cup (D_1 \cap D_2) \in \mathbf{G}$.

(b) Носитель выпуклости — $X$. Поскольку $Y \in \mathbf{H}$, $X = C \cup Y \in \mathbf{G}$ по построению. Минимальным выпуклым элементом выпуклости **G** является $C$, поскольку $\emptyset \in \mathbf{H}$, а по построению $C = C \cup \emptyset$.

(c) $N$-арность. Возьмем произвольное $A = C \cup D \in \mathbf{G}(D \in \mathbf{H})$. Пусть $B \subset A, |B| \leq N$. Разобьем $B = B_C \cup B_D, B_C \subset C, B_D \subset D$. Выпуклая оболочка **G** обладает свойством (следствие из леммы 3, прмененная к **H** и обратному построению):

$$gB = g(B_C \cup B_D) = C \cup hB_D.$$

Поскольку $C$ — минимальное выпуклое, $B_D \subset Y$, то $gB_D = C \cup hB_d$. Так как $|B_D| \leq N, D \in \mathbf{H}$ — $N$-арная, то $hB_D \subset D$. Отсюда

$$gB = C \cup hB_D \subset C \cup D = A.$$

Следовательно, $A$ удовлетворяет определению $N$-арного множества, а потому $A \in \mathbf{G}$.

Обратно, если $A = C \cup D$ удовлетворяет $N$-арному условию, то для любого $B_D \subset D$ с $|B_D| \leq N$ имеем $gB_D \subset A$, откуда $C \cup hB_D \subset C \cup D$, т.е. $hB_D \subset D$. Значит, $D$ удовлетворяет условию $N$-арности в **H**, поскольку **H** — $N$-арная выпуклость, $D \in \mathbf{H}$. Таким образом, **G** в точности состоит из множеств $C \cup D$, $D \in \mathbf{H}$.

(d) Единственность. Если **G** — произвольная $N$-арная выпуклость с $g\emptyset = C$, $\mathbf{H} = \{A \backslash C : A \in \mathbf{G}\}$, то по построению **G** обязано совпадать с $\{C \cup D : D \in \mathbf{H}\}$. Значит, соответствие между **H** и **G** взаимооднозначно.

Таким образом, построенное соответствие устанавливает биекцию между $E_C$ и $\Gamma_0^N(X \backslash C)$. ■

**Предложение.** *Общее число $N$-арных выпуклых структур выражается через биномиальное преобразование над числом заземленных $N$-арных выпуклых структур. А именно*

$$|\Gamma^N(X_n)| = \sum_{k=0}^{n} C_n^k |\Gamma_0^N(X_k)|.$$

**Доказательство.** Классы $E_C$ при различных $C \subset X$ образуют разбиение множества $\Gamma^N(X)$. Мощность каждого класса равна $|\Gamma_0^N(X \backslash C)|$, поэтому

$$|\Gamma^N(X)| = \sum_{C \subset X} |\Gamma_0^N(X \backslash C)|. \quad (7)$$

Случаю $C = \emptyset$ в формуле (7) соответсвует одно слагаемое $|\Gamma_0^N(X_n)|$. Случаю $C = \{x\}, x \in X_n$ в формуле (7) соответвует $n$ таких слагаемых $|\Gamma_0^N(X_{n-1})|$. Аналогично, случаю $|C| = k, 2 \leq k \leq n$, в формуле (7) соответвует $C_n^k$ таких слагаемых $|\Gamma_0^N(X_{n-k})|$ (можно избрать $C_n^k$ различных подмножеств $C \subset X$ мощностью $k$). Таким образом,

$$\sum_{C \subset X} |\Gamma_0^N(X \backslash C)| = \sum_{k=0}^{n} C_n^k |\Gamma_0^N(X_{n-k})|.$$

Проведя замену $k$ на $n-k$ и используя свойство $C_n^k = C_n^{n-k}$, получаем требуемое

$$|\Gamma^N(X)| = \sum_{k=0}^{n} C_n^{n-k} |\Gamma_0^N(X_k)| = \sum_{k=0}^{n} C_n^k |\Gamma_0^N(X_k)|.$$ ■

# 3. ЧИСЛЕННЫЕ ПРИМЕРЫ

На основе разработанного алгоритма и программы был проведён полный перебор $N$-арных выпуклостей на конечных носителях мощности $n = 0,1,\dots,6$ для арностей $N = 0,1,\dots,6$. Результаты численного эксперимента представлены в таблицах 1 и 2 (см. Table 1 и Table 2 в англоязычной версии соответственно). Для $n = 0$ подразумевается единственное пустое множество, а для $n = 1$ — одноэлементный носитель. Для всех исследованных арностей и размерностей выполняется соотношение, установленное в Предложении в разделе 2.

В таблице 3 (см. Table 3 в англоязычной версии) приведены все *бинарные* выпуклости, построенные на носителях с числом элементов 0, 1, 2, 3. Справа от каждой из них после вертикальной черты | указано соотвествующее ей заземление, построенное через исключение минимального выпуклого множества $g\emptyset$ по правилу (1). Пустое множество всюду обозначается символом $\{\}$. В таблице 4 (см. Table 4 в англоязычной версии) приведены все *тернарные,* но *не бинарные* выпуклости, построенные на носителе из 4 элементов.

## СПИСОК ЛИТЕРАТУРЫ